\documentclass[leqno,12pt]{amsart}
\usepackage{latexsym}
\usepackage{mathrsfs}
\usepackage{amsmath}
\usepackage{amssymb}
\usepackage[bookmarks]{hyperref}
\usepackage{pstricks,pst-plot}

\newtheorem{theorem}{Theorem}[section]
\newtheorem{lemma}[theorem]{Lemma}
\newtheorem{corollary}[theorem]{Corollary}
\newtheorem{proposition}[theorem]{Proposition}

\theoremstyle{definition}
\newtheorem{remark}[theorem]{Remark}
\newtheorem{definition}[theorem]{Definition}

\newtheorem{question}[theorem]{Question}

\newcommand{\C}{\mathbb{C}}

\newcommand{\N}{\mathbb{N}}

\newcommand{\Z}{\mathbb{Z}}

\newcommand{\bthm}{\begin{theorem}}
\newcommand{\ethm}{\end{theorem}}
\newcommand{\blem}{\begin{lemma}}
\newcommand{\elem}{\end{lemma}}
\newcommand{\bcor}{\begin{corollary}}
\newcommand{\ecor}{\end{corollary}}
\newcommand{\bprop}{\begin{proposition}}
\newcommand{\eprop}{\end{proposition}}
\newcommand{\bdefn}{\begin{definition}}
\newcommand{\edefn}{\end{definition}}
\newcommand{\bpf}{\begin{proof}}
\newcommand{\epf}{\end{proof}}

\def\vep {\varepsilon}
\def \sm {\setminus}

\def\od{\overline D}

\def\ta {\widetilde A}
\def\tx {\widetilde X}
\def\ti {\widetilde I}

\def\tpi{\widetilde \pi}
\def\ts{\widetilde S}

\def\itemskip {\vskip 3pt plus 2 pt minus 1 pt}
\def\sf{{\mathscr F}}

\def\Asf{A_\sf}
\def\xcf{X\times\C^\sf}

\def\aa{A_\alpha}
\def\ia{I_\alpha}
\def\sa{S_\alpha}
\def\xa{X_\alpha}

\def\piab{\pi_{\alpha,\beta}}
\def\fa{\sf_\alpha}

\def\pimap #1#2{\pi_{#1, #2}}
\def\pistar #1#2{\pi_{#1, #2}^*}
\def\tmap #1#2{T_{#1, #2}}

\newcommand{\ra}{\rightarrow}

\newcommand{\ol}{\overline}

\def\Dk{\Delta_k}
\def\unionDk{\bigcup_{k=1}^\infty \Dk}

\def\singular#1{\exp\left[\displaystyle #1\, \frac{z+1}{z-1}\right]}
\def\singularlam#1{\exp\left[\displaystyle #1\, \frac{z+\lambda}{z-\lambda}\right]}
\def\hull{\mathop{\rm hull}}

\begin{document}
\title[Not strongly regular at a peak point]{A normal uniform algebra that fails\\ to be strongly regular at a peak point}

\author{Alexander J. Izzo}
\address{Department of Mathematics and Statistics, Bowling Green State University, Bowling Green, OH 43403}
\email{aizzo@bgsu.edu}
\thanks{The author was partially supported by NSF Grant DMS-1856010.}


\begin{abstract}
It is shown$\vphantom{\widehat{\widehat{\widehat{\widehat{\widehat{\widehat{\widehat X}}}}}}}$
that there exists a normal uniform algebra, on a compact metrizable space, that fails to be strongly regular at some peak point.  This answers a 31-year-old question of Joel Feinstein.  Our example is $R(K)$ for a certain compact planar set $K$.  Furthermore, it has a totally ordered one-parameter family of closed primary ideals whose hull is a peak point.  General results regarding lifting ideals under Cole root extensions are established.  These results are applied to obtain a normal uniform algebra, on a compact metrizable space, with every point a peak point but again having a totally ordered one-parameter family of closed primary ideals.
\end{abstract}

\maketitle

\vskip -2.54  true in
\centerline{\footnotesize\it Dedicated to Joel Feinstein} 
\vskip 2.54 true in 

%
%
%
%

\section{Introduction}

This paper is devoted to answering questions in the literature regarding strong regularity of uniform algebras and to establishing general results regarding lifting ideals under Cole root extensions.  
Our main goal is to answer the following question raised by Joel Feinstein~\cite[p.~298]{F1} in 1992.
(For definitions of terminology and notation used in this introduction see Section~\ref{notation}.)

\begin{question}\label{question}
 Does there exist a normal uniform algebra $A$, on a compact metrizable space, such that $A$ fails to be strongly regular at some peak point for $A$?  
\end{question}

Part of the interest in this question is that Feinstein showed that it is equivalent to the question: Does there exist a normal uniform algebra $A$, on a compact metrizable space $X$, such that every point of $X$ is a peak point for $A$ but $A$ fails to be strongly regular?
Another closely related question which seems not to be explicitly stated in the literature concerns primary ideals (defined in this context to be those ideals contained in a unique maximal ideal): Does there exist a uniform algebra $A$, normal or not, such that every point of the maximal ideal space of $A$ is a peak point for $A$, but $A$ has a closed primary ideal that is not maximal?
(By a result of Feinstein \cite[Corollary~8]{F2a}, there exist nonnormal uniform algebras such that every point of the maximal ideal space is a peak point.) Question~\ref{question} was 
reiterated by Feinstein with general Banach function algebras in place of uniform algebras in \cite{F2}, and the question also appears, for the specific class of uniform algebras $R(K)$, $K$ a compact planar set, in the forthcoming book of Garth Dales and Ali \" Ulger \cite[Section~3.6]{Dales-Ulger}.  Here, as usual, for $K$ a compact set in the complex plane, $R(K)$ denotes the uniform closure on $K$ of the holomorphic rational functions with poles off $K$.  

We answer all of the above questions affirmatively by establishing the following theorem.

\bthm\label{pre-main-theorem}
There exists a normal uniform algebra $A$, on a compact metrizable space, such that $A$ fails to be strongly regular at some peak point for $A$.  In fact, $A$ can be taken to be $R(K)$ for a certain compact set $K$ in the complex plane.
\ethm

As already mentioned,
Feinstein \cite{F1} showed that the assertion of the first sentence of this theorem is equivalent to the following assertion.

\bcor\label{corollary}
There exists a normal uniform algebra $A$, on a compact metrizable space $X$, such that every point of $X$ is a peak point for $A$ but $A$ is not strongly regular.
\ecor

The inspiration for our proof of Theorem~\ref{pre-main-theorem}, which will give considerably more information than is stated above, comes from the Beurling-Rudin theorem on the closed ideals in the disc algebra \cite{Rudin} (see also \cite[pp.~82--89]{Hoffman}).  Given a compact planar set $K$ contained in the closed unit disc $D$, a point $\lambda$ in $K\cap\partial D$, and a real number $\rho\geq 0$, we will denote by $I_\rho^\lambda$ the closed\vadjust{\kern 5pt} ideal in $R(K)$ generated by the function $(z-\lambda)\singularlam{\rho}$.  In case $\lambda=1$ we will\vadjust{\kern 3pt} write $I_\rho$ in place of $I_\rho^\lambda$.  When $K$ is the closed unit disc, and hence $R(K)$ is the disc algebra, the ideals $I_\rho^\lambda$, $\rho\geq0$, are precisely the closed primary ideals contained in the maximal ideal $M_\lambda$ of functions that vanish at the point $\lambda$.  Furthermore, for $0\leq \rho_1<\rho_2$ there is the strict inclusion $I_{\rho_1}^\lambda\supsetneq I_{\rho_2}^\lambda$.   Our proof of Theorem~\ref{pre-main-theorem} essentially amounts to showing that, taking $\lambda=1$ for instance, McKissick's construction of the first nontrivial normal uniform algebra \cite{McK} (see also \cite[Section~27]{S1}) can be refined so as to preserve this strict inclusion of ideals.  Using results in the author's recent paper \cite{Izzo} we will show that, in addition, the uniform algebra can be chosen in such a way that the point $1$ is the only point where strong regularity fails.
We will thus obtain the following theorem that contains Theorem~\ref{pre-main-theorem}. 
Here, and throughout the paper, we denote the open unit disc in the complex plane by $D$, and given a disc $\Delta$, we denote the radius of $\Delta$ by $r(\Delta)$.

\bthm\label{main-theorem}
For each $r>0$, there exists a sequence of open discs $\{D_k\}_{k=1}^\infty$ such that $\sum_{k=1}^\infty r(D_k)<r$, the point $1$ is in the set $K=\od \sm \bigcup_{k=1}^\infty D_k$, and the following conditions hold:
\begin{enumerate}
\item[(i)] $R(K)$ is normal.
\item[(ii)] $R(K)$ is strongly regular at every point of $K\sm\{1\}$.
\item[(iii)] $R(K)$ is not strongly regular at the point $1$.
\end{enumerate}
Furthermore, the discs $\{D_k\}_{k=1}^\infty$ can be chosen in such a way that $I_{\rho_1} \supsetneq I_{\rho_2}$ for every $0\leq\rho_1<\rho_2$.
\ethm

A modification of the proof of Theorem~\ref{pre-main-theorem} will yield the next result which shows, in particular, that a normal uniform algebra can fail to be strongly regular at an uncountable set of peak points.

\bthm\label{main-theorem-mod}
For each $r>0$, there exists a sequence of open discs $\{D_k\}_{k=1}^\infty$ such that $\sum_{k=1}^\infty r(D_k)<r$ and setting $K=\od \sm \bigcup_{k=1}^\infty D_k$ the following conditions hold:
\begin{enumerate}
\item[(i)] $R(K)$ is normal.
\item[(ii)] $R(K)$ is strongly regular at every point of $K\sm\partial D$.
\item[(iii)] There is a set $\Lambda\subset\partial D$ whose complement in $\partial D$ has one-dimensional Lebesgue measure less than $r$ such that $\Lambda$ is contained in $K$ and at each point of $\Lambda$, the uniform algebra $R(K)$ fails to be strongly regular.
\end{enumerate}
Furthermore, the discs $\{D_k\}_{k=1}^\infty$ can be chosen in such a way that $I_{\rho_1}^\lambda \supsetneq I_{\rho_2}^\lambda$ for every $\lambda\in\Lambda$ and every $0\leq\rho_1<\rho_2$.
\ethm

Note that the uniform algebras in Theorems~\ref{main-theorem} and~\ref{main-theorem-mod}, in spite of failing to be strongly regular, are strongly regular at every \emph{nonpeak point}.
Feinstein and Matthew Heath raised the question of whether there exists a compact planar set $K$ such that $R(K)$ is regular and has no nonzero bounded point derivations, but is not strongly regular \cite[Question~5.8]{FH}.  Each of Theorems~\ref{main-theorem} and~\ref{main-theorem-mod} answers this question affirmatively since for the set $K$ in each of those theorems there are no nonzero bounded point derivations at the points of $K\sm \partial D$ because $R(K)$ is strongly regular at those points, and there are no nonzero point derivations at the points of $K\cap\partial D$ since those points are peak points for $R(K)$.  In \cite{Izzo} the author effectively raised the same question but without the regularity hypothesis, and he promised to give an example answering the question in a future paper.  Thus Theorems~\ref{main-theorem} and~\ref{main-theorem-mod} fulfill that promise.

The next two results show that, as one might expect, Corollary~\ref{corollary} can be strengthened in ways analogous to how Theorems~\ref{main-theorem} and~\ref{main-theorem-mod} strengthen Theorem~\ref{pre-main-theorem}.  

\bthm\label{peak-point}
There exists a normal uniform algebra $B$, on a compact metrizable space $X$, such that every point of $X$ is a peak point for $B$ but there is a point $x_0\in X$ such that there is a one-parameter family $\{H_\rho:0\leq\rho<\infty\}$ of distinct closed primary ideals  contained in the maximal ideal $M_{x_0}$ satisfying $H_{\rho_1}\supsetneq H_{\rho_2}$ for all $0\leq\rho_1<\rho_2$.  Furthermore, $B$ can be taken to have bounded relative units at every point of $X\sm \{x_0\}$.
\ethm

\bthm\label{peak-point-mod}
There exists a normal uniform algebra $B$, on a compact metrizable space $X$, such that every point of $X$ is a peak point for $B$ but there is an uncountable subset $L$ of $X$ such that for every $x\in L$ there is a one-parameter family $\{H_\rho^x: 0\leq\rho<\infty\}$ of distinct closed primary ideals contained in the maximal ideal $M_x$ satisfying $H_{\rho_1}^x\supsetneq H_{\rho_2}^x$ for all $0\leq\rho_1<\rho_2$.  Furthermore, $B$ can be taken to have bounded relative units at every point of $X\sm L$. 
\ethm

Feinstein's proof that Corollary~\ref{corollary} is equivalent to the assertion of the first sentence in Theorem~\ref{pre-main-theorem} used Brian Cole's method of root extensions.
Theorems~\ref{peak-point} and~\ref{peak-point-mod} will be derived from Theorems~\ref{main-theorem} and~\ref{main-theorem-mod} also using Cole's method of root extensions.  However, to do so we will need to prove new results about lifting ideals under root extensions.  There are several examples in the literature in which a certain object lifts under a root extension as a result of adjoining square roots to only a restricted collection of functions.  For instance in~\cite{F1}, \cite{F2}, \cite{F3}, \cite{FH}, \cite{Izzo}, and \cite{ID} square roots are adjoined only to functions that vanish on a given closed set (or a neighborhood of the closed set), and consequently a copy of the given closed set is preserved in the extension.  In~\cite{GI} square roots are adjoined only to functions on which a given bounded point derivation vanishes, with the result that the bounded point derivation lifts to the extended uniform algebra.  
In all these instances, the functions for which square roots are adjoined come from some (proper) closed ideal $I$ of the uniform algebra $A$.  We will prove general results regarding such root extensions.  Very roughly, the results say that in this situation, the quotient Banach algebra $A/I$ is preserved by the extension, and consequently, all the ideals in $A$ that contain $I$ lift under the extension.

In next section we define various terms and notations already used above.  In Section~\ref{preliminaries} we present some known results that we will need.  In Section~\ref{localness-of-strong-reg} we prove that for a normal uniform algebra, strong regularity at a point is a local property.   Although not strictly necesssary for the proofs of Theorems~\ref{main-theorem} and~\ref{main-theorem-mod}, which are presented in Section~\ref{main-proof}, this result greatly simplifies establishing the strong regularity assertion in those theorems.  In Section~\ref{Cole} we present the theorems discussed in the previous paragraph regarding lifting ideals under root extensions.  Theorems~\ref{peak-point} and~\ref{peak-point-mod} are proved in the concluding Section~\ref{peak-point-algebras}.

It is a pleasure to dedicate this paper to Joel Feinstein whose work has been an inspiration to the author and from whom the author has learned some of the notions considered in this paper.

%
%
%
%

\section{Notation and Terminology}\label{notation}

Those readers well versed in uniform algebra concepts may wish to skip or skim this section and refer back to it as needed.

It is to be understood that all sequences, unions, and sums involving an index extend from $1$ to $\infty$; thus for instance $\{D_k\}$ means $\{D_k\}_{k=1}^\infty$, and $\bigcup D_k$ means $\bigcup_{k=1}^\infty D_k$.  
If $f$ is a function whose domain contains a subset $L$, we denote the restriction of $f$ to $L$ by $f|L$, and if $A$ is a collection of such functions, we denote the collection of restrictions of functions in $A$ to $L$ by $A|L$.  
The set of positive integers will be denoted by $\Z_+$.

Throughout the paper all spaces will tacitly be required to be Hausdorff.  Let $X$ be a compact space.  We denote by $C(X)$ the algebra of all continuous complex-valued functions on $X$ equipped with the supremum norm $\|f\|_X=\sup\{|f(x)|: x\in X\}$.  A \emph{uniform algebra} on $X$ is a closed subalgebra of $C(X)$ that contains the constants and separates the points of $X$.  A uniform algebra $A$ on $X$ is said to be 
\itemskip
\begin{enumerate}
\item[(a)] \emph{natural} if the maximal ideal space of $A$ is $X$ (under the usual identification of a point of $X$ with the corresponding multiplicative linear functional),
\itemskip
\item[(b)] \emph{regular on $X$} if for each closed set $K_0$ of $X$ and each point $x$ of $X\sm K_0$, there exists a function $f$ in $A$ such that $f(x)=1$ and $f=0$ on $K_0$,
\itemskip
\item[(c)] \emph{normal on $X$} if for each pair of disjoint closed sets $K_0$ and $K_1$ of $X$, there exists a function $f$ in $A$ such that $f=1$ on $K_1$ and $f=0$ on $K_0$.
\end{enumerate}
The uniform algebra $A$ on $X$ is \emph{regular} or \emph{normal} if $A$ is natural and is regular on $X$ or normal on $X$, respectively.
In fact, every regular uniform algebra is normal \cite[Theorem~27.2]{S1}.
Also, if a uniform algebra $A$ is normal on $X$, then $A$ is necessarily natural \cite[Theorem~27.3]{S1}. 

Let $A$ be a uniform algebra on $X$, and let $x\in X$.  We define the ideals $M_x$ and $J_x$ by
\begin{align} 
M_x &= \{\, f\in A:  f(x)=0\,\}  \nonumber  \\ 
\intertext{and} 
J_x&= \{ \, f\in A: f^{-1}(0)\ \hbox{contains a neighborhood of $x$ in $X$}\}. \nonumber 
\end{align}
More generally, if $E$ is a closed subset of $X$, we define the ideals $M_E$ and $J_E$ by
\begin{align} 
M_E &= \{\, f\in A:  f|E=0\,\}  \nonumber  \\ 
\intertext{and} 
J_E&= \{ \, f\in A: f^{-1}(0)\ \hbox{contains a neighborhood of $E$ in $X$}\}. \nonumber 
\end{align}
When it is necessary to indicate with respect to which algebra the ideals are taken,
we will denote the ideals $J_x$ and $M_x$ in the uniform algebra $A$ by $J_x(A)$ and $M_x(A)$, respectively.

Let $A$ be a natural uniform algebra and let $I$ be an ideal in $A$.   The \emph{hull} of $I$ is the common zero set of the functions in $I$ and is denoted $\hull(I)$.  The ideal $I$ is said to be \emph{local} if $I\supset J\bigl(\hull(I)\bigr)$. The following result is standard \cite[Proposition~4.1.20~(iv)]{D}.

\bthm\label{Garth}
Every ideal in a normal uniform algebra is local.
\ethm

It is also true that every ideal in an $R(K)$ for $K$ a compact planar set is local~\cite[Theorem~1.16]{Izzo}.

The uniform algebra $A$ is \emph{strongly regular at} $x$ if $\ol{J_x}=M_x$, and $A$ is \emph{strongly regular} if $A$ is strongly regular at every point of $X$.  It was shown by Wilken that every strongly regular uniform algebra is normal \cite[Corollary~1]{Wilken2}.

The uniform algebra $A$ has \emph{bounded relative units at $x$} with bound $C\geq 1$ if for each compact subset $K$ of $X\sm \{x\}$, there exists $f\in J_x$ such that $f|K=1$ and $\| f \|_X \leq C$.  If $A$ has bounded relative units at every point of $X$, then $A$ has \emph{bounded relative units}.

The point $x$ is said to be a \emph{peak point} for $A$ if there is a function $f$ in $A$ such that $f(x)=1$ and $|f(y)|<1$ for every $y\in X\sm \{x\}$.  The point $x$ is said to be a \emph{generalized peak point} if for every neighborhood $U$ of $x$ there exists a function $f$ in $A$ such that $f(x)=\|f\|=1$ and $|f(y)|<1$ for every $y\in X\sm U$.  When the space $X$ is metrizable the notions of peak point and generalized peak point coincide.

For $\phi$ a point in the maximal ideal space of the uniform algebra $A$, a \emph{bounded point derivation} on $A$ at $\phi$ is a bounded linear functional $d$ on $A$ satisfying the identity 
\begin{equation*}\label{derivation-eq}
d(fg)=d(f)\phi(g) + \phi(f)d(g)
\end{equation*}
for all $f$ and $g$ in $A$.
It is standard \cite[p.~64]{Browder} that a bounded linear functional $d$ on $A$ is a bounded point derivation at $\phi$ if and only if $d$ annihilates $\ol{M_\phi^2}$ and the constant functions, and hence there exists a bounded point derivation at $\phi$ if and only if $\ol {M_\phi^2}\neq M_\phi$.

%
%
%
%

\section{Preliminaries}\label{preliminaries}

In this section we collect various known results that we will need.  The reader may prefer to skip this section and merely refer back to it when the results are used.

As with McKissick's construction of the first known nontrivial normal uniform algebra \cite{McK}, our proofs of Theorems~\ref{main-theorem} and~\ref{main-theorem-mod} rely on the following lemma.  (Recall that given a disc $\Delta$, we denote the radius of $\Delta$ by $r(\Delta)$.)

\blem\label{McKissick}
Let $\Delta$ be an open disc in the complex plane with center $a$ and radius $r>0$, and let $\vep>0$ be given.  Then there exist a sequence of open discs $\{\Dk\}_{k=1}^\infty$ and a sequence of rational functions $\{f_j\}_{j=1}^\infty$ such that 
\itemskip
\begin{enumerate}
\item[(a)] $\sum_{k=1}^\infty r(\Dk) < \vep$.
\itemskip
\item[(b)] The poles of the $f_j$ lie in $\unionDk$.
\itemskip
\item[(c)] The sequence $\{f_j\}$ converges uniformly on $\C\sm \bigcup_{k=1}^\infty \Dk$ to a function that is identically zero outside $\Delta$ and zero free in $\Delta\sm \unionDk$.
\itemskip
\item[(d)] $\unionDk\subset \{z: r-\vep<|z-a|<r\}$.
\end{enumerate}
\elem

Condition (d) is not part of the lemma as stated by McKissick but is established in the paper \cite{Ko} of Thomas K\" orner where a proof of the lemma simpler than the original one is given.

For the proofs of Theorems~\ref{main-theorem} and~\ref{main-theorem-mod} we will need two recent results of the author \cite[Theorem~1.2 and Lemma~4.2]{Izzo} on strong regularity in $R(K)$ which we state here.

\bthm\label{strongly-regular}
For each $r>0$, there exists a sequence of open discs $\{D_k\}_{k=1}^\infty$ such that $\sum_{k=1}^\infty r(D_k)<r$ and such that setting $K=\ol D\sm \bigcup_{k=1}^\infty D_k$, the uniform algebra $R(K)$ is nontrivial and strongly regular.
\ethm

\blem\label{subsetJ}
Given compact sets $L\subset K\subset \C$ and given a point $x\in L$, if 
$\ol{J_x(R(K))} \supset M_x(R(K))$, then 
$\ol{J_x(R(L))} \supset M_x(R(L))$.
\elem

The following three lemmas will be used to prove condition (v) in Theorem~\ref{Cole-Izzo-thm2}.  The first of these is due to Feinstein and Heath~\cite[Lemma~4.3]{FH}.

\blem\label{bru}
Let $A$ be a uniform algebra on $X$ and $x\in X$.  Suppose that, for each compact subset $E$ of $X\sm \{x\}$, there exists a neighborhood $U$ of $x$ and a function $f\in A$ such that
\begin{enumerate}
\item[(i)] $f|U=1$,
\item[(ii)] $f|E=0$,
\item[(iii)] For each $k\in\N$ there is a function $g\in A$ with $g^{2^k}=f$.
\end{enumerate}
Then $A$ has bounded relative units at $x$.
\elem

The next lemma is part of a result of Feinstein~\cite[Proposition~1.5]{F1}.

\blem\label{bru-implies-peak-point}
Let $A$ be a uniform algebra on a compact space $X$, and let $x\in X$.  If $A$ has bounded relative units at $x$, then $x$ is a generalized peak point for $A$.
\elem

The next lemma, whose elementary proof we omit, is a modification of a lemma of Feinstein~\cite[Lemma~3.5]{F1}.

\blem\label{mod}
Let $A$ be a normal uniform algebra on a compact metrizable space $X$, and let $F$ be a closed subset of $X$.  Then there exists a countable subset $\sf$ of $A$ consisting of functions each vanishing identically on a neighborhood of $F$ such that for each point $x\in X\sm F$, and for each compact subset $E$ of $X\sm\{x\}$, there exists a neighborhood $U$ of $x$, and a function $f\in \sf$ such that $f|U=1$ and $f|E=0$.
\elem

%
%
%
%

\section{Localness of strong regularity}\label{localness-of-strong-reg}

In this section we prove the localness of strong regularity for normal uniform algebras.  
The result will be used in the next section to obtain condition (ii) in Theorems~\ref{main-theorem} and~\ref{main-theorem-mod}.

\bthm\label{localness}
Let $A$ be a normal uniform algebra on a compact space $X$, and let $x_0$ be a point of $X$.  If there exists a closed neighborhood $N$ of $x_0$ in $X$ such that $\ol{A|N}$ is strongly regular at $x_0$, then $A$ is strongly regular at $x_0$.
\ethm

\bpf
Suppose that there exists a closed neighborhood $N$ of $x_0$ in $X$ such that $\ol{A|N}$ is strongly regular at $x_0$.
Fix $f\in A$ satisfying $f(x_0)=0$, and fix $\vep>0$.  We are to show that there exists a function $g\in A$ such that $\|f-g\|_X<\vep$ and $g=0$ on a neighborhood of $x_0$.

Let $L$ be a closed neighborhood of $x_0$ contained in the interior $N^\circ$ of $N$.  By the normality of $A$, there is a function $\varphi\in A$ such that $\varphi=1$ on $L$ and $\varphi=0$ on $X\sm N^\circ$.  

By the strong regularity of $\ol{A|N}$, there is a function $h\in \ol{A|N}$ such that $\|f-h\|_N<\vep/\|\varphi\|_N$ and $h=0$ on some closed neighborhood $M$ of $x_0$.  There is a sequence $(h_n)$ in $A$ such that $h_n|N\ra h$ uniformly on $N$.

For each $n\in\Z_+$, set $g_n=\varphi h_n +(1-\varphi) f$.  Then each $g_n$ is in $A$.  Define a function $g$ on $X$ by
$$
g= 
\begin{cases} \varphi h + (1-\varphi) f &\mbox{on\ }  N \\
f &\mbox{on\ }  X\sm N. 
\end{cases}
$$
Note that on $X\sm N$ we have $g_n=f=g$.  Moreover, $g_n\ra g$ uniformly on $X$.  Thus $g$ is in $A$.  Furthermore, 
\[
\|f-g\|_X=\|f-g\|_N=\|\varphi(f-h)\|_N\leq \|\varphi\|_N\, \|f-h\|_N <\vep.
\]
Finally, $g=0$ on the neighborhood $L\cap M$ of $x_0$.
\epf

%
%
%
%

\section{Proofs of Theorems~\ref{main-theorem} and~\ref{main-theorem-mod}}\label{main-proof}

Recall that given a disc $\Delta$, we denote the radius of $\Delta$ by $r(\Delta)$.  
We will denote the distance from $\Delta$ to the point $1$ by $s(\Delta)$.  Explicitly, $s(\Delta)=\inf\{ |z-1|: z\in \Delta\}$.  We will denote the open disc with center $a_n$ and radius $r_n$ by $D(a_n, r_n)$ and the corresponding closed disc by $\od(a_n, r_n)$.

The following lemma is the key to the proofs of Theorems~\ref{main-theorem} and~\ref{main-theorem-mod}.

\blem
Fix real numbers $0\leq \rho_1<\rho_2$.  Let $\{D_k\}_{k=1}^\infty$ be a sequence of discs such that $\sum_{k=1}^\infty r(D_k) \exp\bigl(2\rho_2/s(D_k)\bigr) <\infty$.  Set $K=\od\sm \bigcup_{k=1}^\infty D_k$.  Then $I_{\rho_1} \supsetneq I_{\rho_2}$.
\elem

\bpf
It is easily shown that the function $(z-1)\singular{\rho_2}$ is in $I_{\rho_1}$, and hence $I_{\rho_1}\supset I_{\rho_2}$.  To prove that the inclusion is strict, we will exhibit a measure on $K$ that annihilates $I_{\rho_2}$ but\vadjust{\kern 1pt} does not annihilate the function $(z-1)\singular{\rho_1}$\vadjust{\kern 6pt}.

For $n\in \Z_+$, let $K_n=\od(0, 1+ \frac{1}{n}) \sm \bigcup_{k=1}^n D_k$.  The boundary of $K_n$ consists of the unions of a finite collection of circular arcs (and possibly some isolated points which can be ignored), and we can define a measure $\mu_n$ on $\partial K_n$ by requiring that for every function $f\in C(K_n)$ we have 
\[
\int\! f\, d\mu_n = \int_{\partial K_n} f(z) \singular{-\rho_2}\, dz.
\]
Let $M=\sum_{k=1}^\infty r(D_k) \exp\bigl(2\rho_2/s(D_k)\bigr) <\infty$.  Then $\|\mu_n\|\leq 2\pi(M+2)$.  Consequently, $\{\mu_n\}_{k=1}^\infty$ has a weak$^*$-accumulation point $\mu$.  Since $K=\bigcap K_n$, the measure $\mu$ is supported on $K$.  If $g$ is a rational function with no poles on $K$, then for large values of $n$, the function $g$ has no poles on $K_n$, and by Cauchy's theorem
\[
\int g(z)(z-1)\singular{\rho_2} \, d\mu_n = \int_{\partial K_n} g(z)(z-1)\, dz = 0.
\]
Thus 
\[
\int g(z)(z-1)\singular{\rho_2} \, d\mu  = 0
\]
for every rational function $g$ with no poles on $K$.  It follows that $\mu$ annihilates $I_{\rho_2}$.

To calculate $\displaystyle\int (z-1)\singular{\rho_1} d\mu$ note that the function\break
$(z-1)\singular {(\rho_1 -\rho_2)}$ has a single isolated singularity at $z=1$, and the residue there is $2(\rho_1-\rho_2)^2 \exp{(\rho_1-\rho_2)}$ since
\begin{align}
(z-&1)\singular{(\rho_1-\rho_2)} = (z-1)\exp\left[ (\rho_1-\rho_2)\left(1+\frac{2}{z-1}\right)\right]\nonumber\\
&=(z-1)\exp\left[ (\rho_1-\rho_2) + \left(\frac{2(\rho_1-\rho_2)}{z-1}\right)\right]\nonumber\\
&=(z-1) \Bigl[\exp{(\rho_1-\rho_2)}\Bigr] \left[1+ \frac{2(\rho_1-\rho_2)}{z-1} + \frac{1}{2}\left(\frac{2(\rho_1-\rho_2)}{z-1}\right)^2 + \cdots \right].\nonumber
\end{align}
Thus by the Residue theorem
\begin{align}
\int(z-1)\singular{\rho_1}\, d\mu_n&=\int_{\partial K_n} (z-1) \singular{(\rho_1-\rho_2)}\, dz 
\nonumber\\
&= 4\pi i(\rho_1-\rho_2)^2 \exp{(\rho_1-\rho_2)}\vphantom{\Biggl[}.
\nonumber
\end{align}
Thus
\[
\int(z-1)\singular{\rho_1} \, d\mu = 4\pi i(\rho_1-\rho_2)^2 \exp{(\rho_1-\rho_2)} \not= 0.
\]
\vskip -8pt 
\epf

\vskip 4pt

\bpf[Proof of Theorem~\ref{main-theorem}]
By the preceding lemma, it suffices to show that discs $\{D_k\}$ can be chosen such that $\sum r(D_k)<r$, such that $\sum r(D_k) \exp\bigl(\nu/s(D_k)\bigr) < \infty$ for every $\nu>0$, and such that conditions (i) and (ii) hold.  We begin by choosing discs such that these conditions are satisfied with the (possible) exception of condition (ii), and then we choose additional discs to achieve condition (ii) in addition.

Choose a sequence $\{\od(a_n,r_n)\}_{n=1}^\infty$ of closed discs such that
\begin{enumerate}
\item[(a)] Each of the discs $\od(1,1/j)$, $j=1, 2, \ldots$, is in $\{\od(a_n, r_n)\}$.
\item[(b)] The discs $\od(1, 1/j)$, $j=1, 2, \ldots$, are the only discs in $\{\od(a_n, r_n)\}$ that contain the point $1$.
\item[(c)] For every $\vep>0$, every point of $\od$ lies in an open disc $D(a_n, r_n)$ with $r_n<\vep$.
\end{enumerate}
Then for each $n\in\Z_+$, the annulus $\{z: r_n/2< |z-a_n|<r_n\}$ is at some positive distance $\delta_n$ from the point $1$.  Set $\vep_n=\min\{2^{-(n+1)}r, 2^{-n}e^{-n/\delta_n}\}$.
For each $n\in\Z_+$, choose discs $\{\Delta_k^n\}_{k=1}^\infty$ in the annulus $\{z: r_n/2<|z-a_n|<r_n\}$ as in Lemma~\ref{McKissick} with $\Delta=D(a_n, r_n)$ and $\vep=\vep_n$.  Then $\sum_{n=1}^\infty\sum_{k=1}^\infty r(\Delta_n^k) < r/2$.  Now let $\nu>0$ be arbitrary.  For each $n\in \Z_+$ we have 
\[
\sum_{k=1}^\infty r(\Delta_k^n) \exp\bigl(\nu/s(\Delta_k^n)\bigr)<\vep_n\exp(\nu/\delta_n).
\]
Thus in particular $\sum_{k=1}^\infty r((\Delta_k^n) \exp(\nu/s\bigl(\Delta_k^n)\bigr)<\infty$.  Furthermore, for all $n>\nu$ we have
\[
\sum_{k=1}^\infty r(\Delta_k^n) \exp\bigl(\nu/s(\Delta_k^n)\bigr) < \vep_n\exp(n/\delta_n)\leq 2^{-n}.
\]
Consequently, $\sum_{n=1}^\infty\sum_{k=1}^\infty r(\Delta_k^n) \exp\bigl(\nu/s(\Delta_k^n)\bigr)<\infty$.

Set $K_1=\od\sm \bigcup_{n=1}^\infty\bigcup_{k=1}^\infty \Delta_k^n$.  Then $R(K_1)$ is regular, and hence normal, for given  a closed set $L\subset K_1$ and a point $x\in K\sm L$, there is some $D(a_n, r_n)$ that contains $x$ and is disjoint from $L$, and hence by the choice of the discs $\{\Delta_k^n\}$ there is a function in $R(K_1)$ that vanishes on $L$ but not at $x$.

To achieve the strong regularity at points different from $1$ we will use Theorem~\ref{strongly-regular} together with Theorem~\ref{localness} on the localness of strong regularity.  Choose a countable collection of open discs $\{B_n\}$ that covers $K_1\sm \{1\}$ such that none of the closed discs $\ol B_n$ contains the point $1$.  Let $\alpha_n$ denote the distance from $\ol B_n$ to the point $1$.  Set $\tilde\vep_n=\min\{2^{-(n+1)}r, 2^{-n}e^{-n/\alpha_n}\}$.  Since in Theorem~\ref{strongly-regular} the open unit disc can, of course, be replaced by any open disc, there exists a sequence of open discs $\{\tilde \Delta_k^n\}$ such that $\sum_{k=1}^\infty r(\tilde \Delta_k^n)<\tilde\vep_n$ and such that setting $K_2^n=\ol B_n\sm \bigcup_{k=1}^\infty \tilde\Delta_k^n$, the uniform algebra $R(K_2^n)$ is strongly regular.  Note that $\sum_{n=1}^\infty\sum_{k=1}^\infty r(\tilde \Delta_k^n)< r/2$ and $\sum_{n=1}^\infty\sum_{k=1}^\infty r(\tilde \Delta_k^n)\exp\bigl(\nu/s(\tilde\Delta_k^n\bigr)<\infty$ for every $\nu>0$ by a computation identical to one done earlier.

Now let $\{D_k\}$ be an enumeration of the collection of discs $\{\Delta_k^n\}_{k,n=1}^\infty \cup \{\tilde\Delta_k^n\}_{k,n=1}^\infty$, and set $K=\od\sm\bigcup D_k$.  Of course $\sum r(D_k) <r$, and $\sum r(D_k) \exp\bigl(\nu/s(D_k)\bigr)<\infty$ for every $\nu>0$.  The uniform algebra $R(K)$ is normal because $K$ is contained in $K_1$ and $R(K_1)$ is normal.  Consider an arbitrary point $x_0\in K\sm\{1\}$, and choose a disc $B_{n_0}$, from the collection $\{B_n\}$, that contains $x_0$.  Then $R(K\cap \ol B_{n_0})$ is strongly regular by Lemma~\ref{subsetJ} because $K\cap \ol B_{n_0}$ is contained in $K_2^{n_0}$ and $R(K_2^{n_0})$ is strongly regular.  Furthermore, $\ol{R(K)|(K\cap\ol B_{n_0})}=R(K\cap \ol B_{n_0})$.  Thus Theorem~\ref{localness} shows that $R(K)$ is strongly regular at $x_0$.
\epf

\bpf[Proof of Theorem~\ref{main-theorem-mod}]
Since the proof is similar to the proof of Theorem~\ref{main-theorem}, we merely indicate the modifications needed.  For the discs $\{\od(a_n, r_n)\}$, we discard conditions (a) and (b) and instead require that each disc $\od(a_n, r_n)$ is either centered at a point of $\partial D$ or else is contained in $D$, and we retain condition (c).  Let $\Gamma$ denote the set of $n$ such that $a_n$ is in $\partial D$.  For each $n\in \Gamma$, choose a number $\gamma_n>0$ in such a way that the intersection with $\partial D$ of the union of the annuli $\{z: r_n-\gamma_n < |z-a_n| < r_n+\gamma_n\}$ has one-dimensional Lebesgue measure less than $r$.  Let $\Lambda$ be the complement of that union in $\partial D$.
Choose, for each $n\in\Z_+$, a positive number $r_n'$ such that $r_n-\gamma_n<r_n'<r_n$.  Then the annulus $\{z:r_n'<|z-a_n|<r_n\}$ is at a positive distance $\delta_n$ from $\Lambda$.  
To establish everything except condition (ii), we choose discs $\{\Delta_k^n\}$ in the annulus $\{z:r_n'<|z-a_n|<r_n\}$ as in Lemma~\ref{McKissick} arguing as in the proof of Theorem~\ref{main-theorem}, but with distance $s(\Delta_k^n)$ to $1$ replaced by distance to $\Lambda$.  

To get condition (ii), we argue essentially as in the proof Theorem~\ref{main-theorem} except that for the collection $\{B_n\}$ we take the collection $\{D(0, 1- \frac{1}{n}): n=2, 3, \ldots\}$, and we again replace distance to $1$ by distance to $\Lambda$.  
\epf


%
%
%
%

\section{Root extensions and ideals}\label{Cole}

In this section we prove results about systems of root extensions and ideals which we will use in the next section to prove Theorems~\ref{peak-point} and~\ref{peak-point-mod}.  We present the results in greater generality than we will need because we believe they are of interest in their own right and are likely to have further applications.

Cole's method of root extensions involves an iterative process \cite{Cole}. We begin by discussing a single step of the iteration.

Let $A$ be a uniform algebra on a compact space $X$, and let $\sf$ be a (nonempty) subset of $A$.  Endow $\C^\sf$ with the product topology.  Let $p_1:X\times \C^\sf\ra X$ and $p_f:\xcf\ra\C$ denote the projections given by $p_1(x, (z_g)_{g\in \sf})=x$ and $p_f(x, (z_g)_{g\in\sf})=z_f$.  Define $X_\sf\subset \xcf$ by
\[
X_\sf = \{\, y\in\xcf:\bigl(p_f(y)\bigr)^2=f\bigl(p_1(y)\bigr) \ \hbox{for all $f\in\sf$}\,\},
\]
and let $A_\sf$ be the uniform algebra on $X_\sf$ generated by the set of functions $\{\, f\circ p_1: f\in A\} \cup \{\, p_f: f\in\sf\}$.  On $X_\sf$ we have $p_f^2=f\circ p_1$ for every $f\in \sf$.  Set $\pi=p_1|X_\sf$, and note that $\pi$ is surjective.  There is an isometric embedding $\pi^*:A\ra A_\sf$ given by $\pi^*(f)=f\circ \pi$.

We call the uniform algebra $A_\sf$ or the pair $(A_\sf, X_\sf)$, the $\sf$-extension of $A$, and we call $\pi$ the associated surjection.  Note that if $X$ is metrizable and $\sf$ is countable, then $X_\sf$ is metrizable also.  
Given $x\in X$, if $\sf$ is contained in $M_x$, then the set 
$\pi^{-1}(x)$ consists of a single point.  

There is an operator $S:A_\sf\ra \pi^*(A)$ given by integrating over the fibers of $\pi$ using the measure on each fiber that is invariant under the obvious action of $(\Z/2)^\sf$ on each fiber.  See \cite{Cole} or \cite[pp.~194--195]{S1} for details.  
Rather than working with $S$, we will use the operator 
$T: A_\sf\ra A$ obtained from $S$ by identifying $\pi^*(A)$ with $A$.  
The following properties of $T$ are almost obvious.

\blem\label{key-trick}
\phantom{.}
\begin{enumerate}
\item[(i)] $\|T\|=1$.
\item[(ii)] $T\circ \pi^*$ is the identity.
\item[(iii)] Given distinct functions $f_1,\ldots, f_r\in \sf$ and a function $f\in A$,
$$T\bigl(\pi^*(f) p_{f_1}\cdots p_{f_r}\bigr)=0.$$
\end{enumerate}
\elem

One can iterate the above extension process.  This leads to the notion of a system of root extensions which we next define.

Henceforth, $\tau$ will be a fixed infinite ordinal.  
A \emph{system of root extensions} is a triple of indexed sets $\bigl(\{\aa\}, \{\xa\}, \{\piab\}\bigr)$ $(0\leq \alpha\leq\beta\leq\tau)$ (denoted for brevity by $\{\aa\}_{0\leq\alpha\leq\tau}$) where
each $\xa$ is a compact space, each $\aa$ is a uniform algebra on $\xa$, and each $\piab$ is a continuous surjective map $\piab:X_\beta\ra \xa$ such that the following conditions hold:
\begin{enumerate}
\item[(i)] The equation $\pi^*_{\alpha,\beta}(f)=f\circ \piab$ defines a homomorphism of $\aa$ into $A_\beta$.
\item[(ii)] For $\alpha\leq\beta\leq\gamma$, $\pi_{\alpha,\beta}\circ\pi_{\beta,\gamma}=\pi_{\alpha,\gamma}$, and $\pi_{\alpha,\alpha}$ is the identity on $\xa$.
\item[(iii)] For $\alpha<\tau$, there is a subset $\fa$ of $\aa$ such that $A_{\alpha+1}$ is the $\fa$-extension of $\aa$ and $\pi_{\alpha,\alpha+1}$ is the associated surjection.
\item[(iv)] For $\gamma$ a limit ordinal, $X_\gamma$ is the inverse limit of the inverse system $\{\xa, \piab\}_{\alpha\leq\beta<\gamma}$, the maps $\pi_{\alpha,\gamma}:X_\gamma\ra\xa$ are those associated with the inverse limit, and $A_\gamma$ is the closure in $C(X_\gamma)$ of $\bigcup\limits_{\alpha<\gamma} \pi^*_{\alpha,\gamma}(\aa)$.
\end{enumerate}

The existence of systems of root extensions is of course proved by transfinite induction.  A choice of the subsets $\sf_\alpha$ uniquely determines a system of root extensions.

\begin{remark}\label{remark}
It follows trivially from conditions (i) and (ii) that for $\alpha\leq\beta\leq\gamma$, $\pistar\beta\gamma\circ\pistar\alpha\beta=\pistar\alpha\gamma$, and $\pistar\alpha\alpha$ is the identity on $\aa$.
\end{remark}

Given a uniform algebra $A$ on $X$, a uniform algebra $\ta$ on $\tx$, and a surjective continuous map $\tpi:\tx\ra X$, we will say that $\ta$ and $\tpi$ are obtained from $A$ by a system of root extensions if there exists a system of root extensions $\bigl(\{\aa\}, \{\xa\}, \{\piab\}\bigr)$ $(0\leq \alpha\leq\beta\leq\tau)$ with $A_0=A$, $A_\tau=\ta$, and $\pimap 0\tau=\tpi$.

The following is \cite[Corollary~2.9]{F1} of Feinstein.

\blem\label{preserve-normal2}
Given a system of root extensions $\{A_\alpha\}_{0\leq \alpha\leq \tau}$, if $A_0$ is normal, then $A_\alpha$ is normal for all $\alpha$.
\elem

For a system of root extensions $\{\aa\}_{0\leq\alpha\leq\tau}$, Cole introduced certain surjective linear operators $T_\beta:A_\beta\ra A_0$.  It will be helpful for us to introduce, more generally, operators $\tmap\alpha\beta: A_\beta\ra \aa$ for every $\alpha\leq\beta$.

\blem\label{T-op}
Given a system of root extensions $\{A_\alpha\}_{0\leq \alpha\leq \tau}$ there exists a system of surjective linear operators $\{\tmap\alpha\beta: A_\beta\ra \aa\}_{0\leq\alpha\leq\beta\leq\tau}$ such that for all $0\leq\alpha\leq\beta\leq\gamma\leq\tau$
\begin{enumerate}
\item[(a)] $\tmap\alpha\alpha$ is the identity operator on $\aa$.
\item[(b)] $\|\tmap\alpha\beta\|=1$.
\item[(c)] $\tmap\alpha\beta\circ\tmap\beta\gamma=\tmap\alpha\gamma$.
\item[(d)] $\tmap\alpha\beta\circ\pistar\alpha\beta$ is the identity on $\aa$.
\item[(e)] $\tmap\alpha\gamma\circ\pistar\beta\gamma=\tmap\alpha\beta$.
\item[(f)]$\tmap\beta\gamma\circ\pistar\alpha\gamma=\pistar\alpha\beta$.
\end{enumerate}
\elem

\bpf
First note that condition (d) is an immediate consequence of conditions (a) and (e), and note that condition (f) is an immediate consequence of condition (d) and Remark~\ref{remark}.  Thus is suffices to show that the operators $\{\tmap\alpha\beta\}$ can be chosen to satisfy conditions (a), (b), (c), and (e).  We will apply transfinite induction on $\beta$ with $\alpha$ fixed to obtain operators $\{\tmap\alpha\beta\}$ satisfying conditions (a), (b), and (e), and then observe that these operators satisfy condition (c) also.

The operator $\tmap\alpha\alpha$ is specified.  Consider $\alpha\leq\beta\leq\tau$, and assume as the induction hypothesis that operators $\tmap\alpha\delta$ have been defined for all $\alpha\leq\delta<\beta$ such that conditions (a), (b), and (e) hold.  If $\beta=\delta+1$ for some $\delta$, then $A_\beta$ is the $\sf_\beta$-extension of $A_\delta$.  Let $T:A_\beta\ra\ A_\delta$ be the operator discussed in the paragraph immediately preceding Lemma~\ref{key-trick}, set $\tmap\alpha\beta=\tmap\alpha\delta\circ T$, and verify that conditions (a), (b), and (e) continue to hold.  If $\beta$ is a limit ordinal, define an operator 
$\widetilde T_{\alpha,\beta}$ on the dense subspace $\bigcup_{\delta<\beta} \pistar\delta\beta(A_\delta)$ of $A_\beta$ by 
\[
\widetilde T_{\alpha,\beta}(\pistar\delta\beta f)=\tmap\alpha\delta f \ \hbox{for} \ f\in A_\delta.
\]
Condition (e) insures that $\widetilde T_{\alpha,\beta}$ is well defined, i.e., if $f_1\in A_{\delta_1}$ and $f_2\in A_{\delta_2}$ satisfy $\pistar{\delta_1}\beta f_1 = \pistar{\delta_2}\beta f_2$, then $\tmap{\alpha}{\delta_1} f_1 = \tmap\alpha{\delta_2} f_2$.  Furthermore, $\|\widetilde T_{\alpha,\beta}\|=1$, so $\widetilde T_{\alpha,\beta}$ has a unique continuous extension to an operator on $A_\beta$ which we declare to be $\tmap\alpha\beta$.  Conditions (a), (b), and (e) continue to hold.  Thus the existence of operators $\{\tmap\alpha\beta\}_{0\leq\alpha\leq\beta\leq\tau}$ satisfying conditions (a), (b), and (e) is established.

To verify that the operators we have defined satisfy condition (c),  fix $\alpha$ and $\beta$, and apply transfinite induction on $\gamma$.
\epf

\blem\label{one-pt-fiber}
If $\pi^{-1}_{\alpha,\beta}(x)$ consists of a single point, then $(\tmap\alpha\beta f)(x)=f(\pi^{-1}_{\alpha,\beta}(x))$ for each function $f\in C(X_\beta)$.
\elem

\bpf
For fixed $\alpha$, apply transfinite induction on $\beta$.
\epf

We will need the following functional analysis lemma whose elementary proof we omit.

\blem\label{elementary}
Let $X$ be a Banach space, let $Y$ be a closed subspace of $X$, let $I$ be a closed subspace of $Y$, and let $S:X\ra X$ be a norm $1$ projection of $X$ onto $Y$.  Then the map $\tilde S:X/T^{-1}(I)\ra X/I$ induced by $S$ is an isometry.
\elem

The next lemma is the key to the proofs of our results on root extensions and ideals.  Its proof is similar to the proof of \cite[Lemma~4.1]{GI} which is essentially the special case in which the closed ideal $I$ arises from a bounded point derivation.  In the lemma, $T:A_\sf\ra A$ is the operator in Lemma~\ref{key-trick}.

\blem\label{key-Cole}
Let $A$ be a uniform algebra on a compact space $X$, let $I$ be a closed ideal in $A$, and let  
$\phi:A\ra A/I$ denote the quotient map.  Let
$\sf$ be a subset of $I$.  
Then in the $\sf$-extension $A_\sf$ of $A$ the set $I_\sf=T^{-1}(I)$ is a closed ideal, and the map $\phi\circ T:A_\sf\ra A/I$ is a Banach algebra homomorphism that induces an isometric Banach algebra isomorphism of $A_\sf/I_\sf$ onto $A/I$.
\elem

\bpf
For notational convenience set $\Phi=\phi\circ T$.  Then $\Phi$ is a linear map with kernel $I_\sf$.  Therefore, if $\Phi$ is multiplicative, i.e., satisfies
\begin{equation}\label{multiplicative}
\Phi(fg)=\Phi(f)\Phi(g)
\end{equation}
for all $f,g\in A_\sf$, then $\Phi$ is a Banach algebra homomorphism, $I_\sf$ is an ideal in $A_\sf$, and identifying $\pi^*(A)$ with $A_\sf$ and applying Lemma~\ref{elementary} shows that the induced Banach algebra isomorphism of $A_\sf/I_\sf$ onto $A/I$ is isometric.
Thus it suffices to show that $\Phi$ satisfies equation (\ref{multiplicative}) for all $f,g\in A_\sf$.
Moreover, it is enough to verify equation~(\ref{multiplicative}) for $f$ and $g$ belonging to
the dense subalgebra $H$ of $\Asf$ that is algebraically generated by $\pi^*(A) \cup \{\, p_f:f\in \sf\}\, $.  Functions $f$ and $g$ in $H$ can be expressed in the form 
\[
f = \pi^*(f_0)+\sum_{u=1}^s \pi^*(f_u) F_u \qquad\hbox{and} \qquad g = \pi^*(g_0)+\sum_{v=1}^t \pi^*(g_v) G_v
\]
where $f_0, f_1,\ldots,f_s, g_0,g_1,\ldots, g_t\in A$ and each $F_u$ and each $G_v$ is a nonempty product of distinct functions of the form $p_f$ for $f\in \sf$.

By Lemma~\ref{key-trick}, $Tf=f_0$ and $Tg=g_0$, so
\[
\Phi(f)=(\phi\circ T)(f)=\phi(f_0) \qquad\hbox{and} \qquad \Phi(g)=(\phi\circ T)(g)=\phi(g_0).
\]
Since $\phi(f_0g_0)=\phi(f_0)\phi(g_0)$,
the proof will be complete once we show that $\Phi(fg)=\phi(f_0g_0)$.

View $fg$ as a sum of four terms:
\begin{align}
fg =\pi^*(f_0g_0) &+ \biggl(\sum_{u=1}^s \pi^*(f_ug_0) F_u\biggr)
+ \biggl(\sum_{v=1}^t \pi^*(f_0g_v) G_v\biggr)\nonumber\\
&+ \biggl(\sum_{u=1}^s\sum_{v=1}^t \pi^*(f_ug_v)F_uG_v\biggr).\nonumber
\end{align}
By Lemma~\ref{key-trick},
\begin{equation}\label{part1}
T\bigl(\pi^*(f_0g_0)\bigr)=f_0g_0
\end{equation}
\begin{equation}\label{part2}
T\biggl(\sum_{u=1}^s \pi^*(f_ug_0) F_u\biggr)=0
\end{equation}
\begin{equation}\label{part3}
T\biggl(\sum_{v=1}^t \pi^*(f_0g_v) G_v\biggr)=0.
\end{equation}
Now for fixed $u$ and $v$, consider $T\bigl(\pi^*(f_ug_v)F_uG_v\bigr)$.
We have $F_u=p_{f_1}\cdots p_{f_a}$ and $G_v=p_{g_1}\cdots p_{g_b}$ where $f_1,\ldots,f_a$ are distinct elements of $\sf$ and $g_1,\ldots, g_b$ are also distinct elements of $\sf$.  Note that each of the sets $\{ f_1,\ldots, f_a\}$ and $\{g_1,\ldots, g_b\}$ is necessarily nonempty. If $\{ f_1,\ldots, f_a\}=\{g_1,\ldots, g_b\}$, then $F_uG_v=p_{f_1}^2\cdots p_{f_a}^2=\pi^*(f_1\cdots f_a)$, and hence
\[
(\phi\circ T)\bigl(\pi^*(f_ug_v)F_uG_v\bigr) = (\phi\circ T)\bigl(\pi^*(f_ug_v f_1\cdots f_a)\bigr) =  \phi(f_ug_v f_1\cdots f_a);
\]
the last quantity above is zero because $f_1,\ldots, f_a$ belong to the ideal $I$.  If instead $\{f_1,\ldots, f_a\}\neq \{g_1,\ldots, g_b\}$, then $F_uG_v$ can be expressed as the product of a possibly empty set of elements of $\pi^*(A)$ and a \emph{nonempty} set of functions $p_{h_1},\ldots, p_{h_c}$ with $h_1,\ldots, h_c\in \{f_1,\ldots, f_a, g_1,\ldots, g_b\}$; consequently,
$T\bigl(\pi^*(f_ug_v)F_uG_v\bigr)=0$ by Lemma~\ref{key-trick}(iii).  We conclude that 
\begin{equation}\label{part4}
(\phi\circ T)\biggl(\sum_{u=1}^s\sum_{v=1}^t \pi^*(f_ug_v)F_uG_v\biggr)=0.
\end{equation}
Collectively, equations (\ref{part1})--(\ref{part4}) yield that
\[
\Phi(fg) =  (\phi\circ T)(fg) = \phi(f_0g_0),
\]
as desired.
\epf

Finally we come to the theorems of this section.

\bthm\label{Cole-Izzo-thm1}
Let $\bigl(\{\aa\}, \{\xa\}, \{\piab\}\bigr)$ $(0\leq \alpha\leq\beta\leq\tau)$ be a system of root extensions.
Let $I_0$ be a closed ideal in $A_0$, and set $S_0=\hull(I_0)$.
For every $0\leq \alpha\leq\tau$, set $\ia=T_{0, \alpha}^{-1}(I_0)$ and $\sa=\pi_{0, \alpha}^{-1}(S_0)$.  Suppose that $\ia\supset\fa$ for every $0\leq\alpha<\tau$.
Then for every $0\leq\alpha\leq\tau$
\begin{enumerate}
\item $\pimap 0\alpha$ takes $\sa$ homeomorphically onto $S_0$.
\item $\pistar 0\alpha$ induces in the obvious way an isometric isomorphism of $A_0|S_0$ onto $\aa|\sa$.
\item $\ia$ is a closed ideal in $\aa$ such that $\hull(\ia)=\sa$, such that $\ia\cap\pistar 0\alpha(A)=\pistar 0\alpha(I_0)$, and such that $\aa/\ia$ is isometrically isomorphic as a Banach algebra to $A_0/I_0$.  Consequently, there is an order-preserving bijective correspondence between the closed ideals of $\aa$ containing $\ia$ and the closed ideals of $A_0$ containing $I_0$.
\end{enumerate}
\ethm

\bpf
By a simple transfinite induction, one shows simultaneously that the hypotheses imply that every function in $\ia$, and hence every function in $\sf_\alpha$, is zero on $\sa$, and that 
$\pimap 0\alpha$ takes $\sa$ one-to-one onto $S_0$.  Since $\pimap 0\alpha$ is continuous and $\sa$ is compact, condition (i) follows.

Given $f\in A_0$, the restriction of $\pistar 0\alpha(f)=f\circ \pimap 0\alpha$ to $\pi^{-1}_{0,\alpha}(S)=S_\alpha$ depends only on the restriction of $f$ to $S_0$, so $\pistar 0\alpha$ induces a map of $A_0|S_0$ into $\aa|\sa$ that is obviously isometric.  Given $g\in \aa$, Lemma~\ref{one-pt-fiber} shows that $\pistar 0\alpha(\tmap 0\alpha g)|\sa=g$, so the map is onto.  Thus condition (ii) holds.

We have already noted that every function in $\ia$ is zero on $\sa$.  Since for every point $x$ of $\xa\sm\sa$ there is a function $f\in I_0$ such that $\pistar0\alpha(f)$, a function in $\ia$, is nonzero at $x$, this gives that $\hull(\ia)=\sa$.

The equality $\ia\cap\pistar 0\alpha(A_0)=\pistar 0\alpha(I_0)$ follows immediately from Lemma~\ref{T-op}~(d).

It remains to show, for each $\alpha$, that $\ia$ is an ideal in $\aa$ and that $\aa/\ia$ and $A_0/I_0$ are isometrically isomorphic as Banach algebras.  Let $\phi:A_0\ra A_0/I_0$ denote the quotient map.  Assume for the moment that the map $\phi\circ \tmap 0\alpha: \aa\ra A_0/I_0$ is a Banach algebra homomorphism for every $\alpha$.  Then since $\ia=\ker (\phi\circ\tmap 0\alpha)$, it follows that $\ia$ is an ideal in $\aa$, and that the induced map $\aa/\ia\ra A_0/I_0$ is a Banach algebra isomorphism.  Identifying $A_0$ with the subspace $\pistar 0\alpha(A_0)$ of $\aa$ and applying Lemma~\ref{elementary} shows that this isomorphism is an isometry.
Thus to complete the proof it suffices to show that the map $\phi\circ \tmap 0\alpha: \aa\ra A_0/I_0$ is indeed a Banach algebra homomorphism for every $\alpha$.

\clubpenalty=0
\widowpenalty=0
We apply transfinite induction.  Consider $0\leq \beta\leq\tau$, and assume as the induction hypothesis that $\phi\circ\tmap0\alpha$ is a Banach algebra homomorphism for every $\alpha<\beta$.  In case $\beta=0$, nothing needs to be proved.  If $\beta$ is a limit ordinal, then it is immediate from the induction hypothesis that the restriction of $\phi\circ\tmap0\beta$ to the dense subset $\bigcup_{\alpha<\beta} \pistar\alpha\beta(\aa)$ of $A_\beta$ is an algebra homomorphism, and hence, $\phi\circ\tmap0\beta$ is a Banach algebra homomorphism by continuity.  Now suppose instead that $\beta=\gamma+1$ for some $\gamma$.  The map $\phi\circ\tmap0\gamma$ is a Banach algebra homomorphism by the induction hypothesis.  Consequently, $I_\gamma$ is a closed ideal in $A_\gamma$.  Let $\phi_\gamma:A_\gamma\ra A_\gamma/I_\gamma$ denote the quotient map, and let $\iota_\gamma:A_\gamma/I_\gamma\ra A_0/I_0$ denote the Banach algebra isomorphism induced by $\phi\circ\tmap0\gamma$.  By Lemma~\ref{key-Cole}, the map $\phi_\gamma\circ\tmap\gamma{\gamma+1}: A_{\gamma+1}\ra A_\gamma/ I\gamma$ is a Banach algebra homomorphism.  Now consider the commutative diagram 
\vskip 32pt
\[ 
\begin{array}{ccccc}
A_{\gamma+1}\phantom{\Big[} & \stackrel{\textstyle\tmap\gamma{\gamma+1}}{\relbar\joinrel\relbar\joinrel\relbar\joinrel\longrightarrow} &A_\gamma 
& \stackrel{\textstyle\tmap0\gamma}{\relbar\joinrel\relbar\joinrel\relbar\joinrel\longrightarrow} &A_0\\ 
\bigg\downarrow && 
\bigg\downarrow\rlap{${\phi_\gamma}\,\,\,$}&& 
\bigg\downarrow\rlap{${\phi}\,\,\,$}
\\
A_{\gamma+1}/I_{\gamma+1} & \stackrel{}{\relbar\joinrel\relbar\joinrel\relbar\joinrel\longrightarrow} &A_\gamma/I_\gamma
& \stackrel{\textstyle\iota_\gamma}{\relbar\joinrel\relbar\joinrel\relbar\joinrel\longrightarrow} 
&A_0/I_0.
\\
\end{array} 
\]
%
\psset{linewidth=0.45pt}
\psset{xunit=.76cm}
\psset{yunit=.76cm}
\psset{runit=.76cm}
\vphantom{.}
\vskip 0.035 in
\psline[arrowsize=4pt]{->}(5.2,3.3)(7.1,2.1)  
\psline[arrowsize=4pt]{->}(9.5,3.3)(11.4,2.1)  
\pscurve[arrowsize=4pt]{->}(4.0,4.6)(6.0,5.2)(8.0,5.4)(10.0,5.2)(12.0,4.6)
\pscurve[arrowsize=4pt]{->}(4.0,0.8)(6.0,0.2)(8.0,0.0)(10.0,0.2)(12.0,0.8)
\vskip -1.94in
\hskip 2.38 in$\tmap0{\gamma+1}$
\vskip 2.0 in
\psset{linewidth=0.8pt}
\noindent
Observe that the map $\phi\circ\tmap0{\gamma+1}: A_{\gamma+1}\ra A_0/I_0$ coincides with the composition of the Banach algebra homomorphisms $\phi_\gamma\circ\tmap\gamma{\gamma+1}$ and $\iota_\gamma$ and hence is itself a Banach algebra homomorphism, as desired.
\epf

By suitable choice of system of root extensions we will obtain the following as a corollary.

\clubpenalty=150
\widowpenalty=150

\bthm\label{Cole-Izzo-thm2}
Let $A$ be a uniform algebra on a compact space $X$, let $I$ be an ideal in $A$, and set $S=\hull(I)$.
Then there exists a uniform algebra $\ta$ on a compact space $\tx$ and a surjective continuous map $\tpi:\tx\ra X$, obtained from $A$ by a system of root extensions, and there exists an ideal $\ti$ in $\ta$ such that setting $\ts=\tpi^{-1}(S)$ conditions (i)--(iii) of Theorem~\ref{Cole-Izzo-thm1} hold with $\ta$, $\tx$, $\ti$, $\ts$, and $\tpi$ in place of $\aa$, $\xa$, $\ia$, $\sa$, and $\pimap 0\alpha$, respectively, and such that furthermore,
\begin{enumerate}
\setcounter{enumi}{3}
\item Every function in $\ti$ has a square root in $\ti$.
\item If $A$ is normal, then $\ta$ is normal and has bounded relative units at every point of $\tx\sm\ts$, and hence every point of $\tx\sm\ts$ is a generalized peak point for $\ta$.
\end{enumerate}
If $X$ is metrizable, then, in addition, we can take $\tx$ to be metrizable provided we replace condition (iv) by 
\begin{enumerate}
\item[(iv$'$)] there is a dense subset $\sf$ of $\ti$ such that every function in $\sf$ has a square root in $\sf$.
\end{enumerate}
\ethm

\bpf
Using Theorem~\ref{Cole-Izzo-thm1} and transfinite induction, it is easily shown that there is a system of root extensions satisfying the conditions of Theorem~\ref{Cole-Izzo-thm1} with $\tau=\Omega$ (the first uncountable ordinal) and $\sf_\alpha=\ia$ for every $0\leq \alpha<\Omega$.  Set $\ta=A_\Omega$, $\tx=X_\Omega$, $\ti=I_\Omega$, etc.  Then conditions (i)--(iii) hold with $\ta$, $\tx$, $\ti$, $\ldots$  in place of $\aa$, $\xa$, $\ia$, $\ldots$, respectively, by Theorem~\ref{Cole-Izzo-thm1}.

Given $f\in\ti=I_\Omega$, there is some $\alpha<\Omega$ and some $g\in \aa$ such that $f=\pistar \alpha\Omega g$.  By construction and Lemma~\ref{key-trick}~(iii), $\pistar \alpha{\alpha+1} g=h^2$ for some $h\in A_{\alpha+1}$ such that 
\begin{equation}
\tmap \alpha{\alpha+1} h=0.
\end{equation}\label{simple}
Now
\[
\left(\pistar{\alpha+1}\Omega h\right)^2 =\pistar{\alpha+1}\Omega h^2 =\pistar{\alpha+1}\Omega \pistar\alpha{\alpha+1} g = \pistar\alpha\Omega g = f.
\]
Furthermore, by Lemma~\ref{T-op} and equation~(\ref{simple})
\[
\tmap 0\Omega\left(\pistar{\alpha+1}\Omega h\right) = \left(\tmap 0{\alpha+1} \circ \tmap{\alpha+1}\Omega\right) \left(\pistar{\alpha+1}\Omega h\right) = \tmap 0{\alpha+1} h =
\tmap0\alpha\circ\tmap\alpha{\alpha+1} h= 0,
\]
so $\pistar {\alpha+1}\Omega h$ is in $\ti$.  Thus every function in $\ti$ has a square root in $\ti$.

Now suppose that $A$ is normal.  Then by Lemma~\ref{preserve-normal2}, $\ta$ is normal.  Consequently, given a point $\tilde x\in \tx\sm\ts$ and a compact subset $E$ of $\tx\sm\{\tilde x\}$, Theorem~\ref{Garth} insures that there is a function in $\ti$ that is one on a neighborhood of $\tilde x$ and zero on $E$.  Therefore, $\ta$ has bounded relative units at $\tilde x$ by Lemma~\ref{bru}.  The final assertion of condition (v) follows by Lemma~\ref{bru-implies-peak-point}.

All that remains is to prove the last sentence of the theorem.  From now on suppose that $X$ is metrizable.  Using Theorem~ \ref{Cole-Izzo-thm1} and transfinite induction, it is easily shown that there is a system of root extensions satisfying the conditions of Theorem~\ref{Cole-Izzo-thm1} with $\tau=\omega$ (the first infinite ordinal) and such that, for every $0\leq \alpha<\omega$, the collection $\sf_\alpha$ a countable dense subset of $\ia$ such that for every function $f\in\sf_\alpha$ the function $\pistar \alpha {\alpha+1} f$ is the square of a function in $\sf_{\alpha+1}$.  Furthermore, if $A$ is normal, then Lemma~\ref{mod} and Theorem~\ref{Garth} show that we can, and therefore we shall, choose $\sf_\alpha$ such that,
setting $\sa=\pi_{0,\alpha}^{-1}(S)$, we have that for each point $x\in \xa\sm\sa$, and for each compact subset $E$ of $\xa\sm\{x\}$, there exists a neighborhood $U$ of $x$, and a function $f\in \sf_\alpha$ such that $f|U=1$ and $f|E=0$.  
Set $\ta=A_\omega$, $\tx=X_\omega$, $\ti=I_\omega$, etc.  
Then $\tx$ is metrizable.
Furthermore, conditions (i)--(iii) hold with $\ta$, $\tx$, $\ti$, $\ldots$  in place of $\aa$, $\xa$, $\ia$, 
$\ldots$, respectively, by Theorem~\ref{Cole-Izzo-thm1}.  

We will establish condition (iv$'$) with 
$\sf=\bigcup_{0\leq\alpha<\omega} \pistar \alpha\omega (\sf_\alpha)$.  
First we show that every function in $\sf$ has a square root in $\sf$.  Given $f\in\sf$, there is some $\alpha<\omega$ and some $g\in \sf_\alpha$ such that $f=\pistar \alpha \omega g$.  By construction, $\pistar \alpha{\alpha+1} g=h^2$ for some $h\in\sf_{\alpha+1}$.  Then $\pistar {\alpha+1}\omega h$ is in $\sf$ and by Remark~\ref{remark}
\[
\left(\pistar{\alpha+1}\omega h\right)^2 =\pistar{\alpha+1}\omega h^2 =\pistar{\alpha+1}\omega \pistar\alpha{\alpha+1} g = \pistar\alpha\omega g = f,
\]
so $f$ has a square root in $\sf$.

Next we show that $\sf$ is contained in $\ti$.  Let $f$, $\alpha$, and $g$ be as in the previous paragraph.  Then by Lemma~\ref{T-op}
\[
\tmap 0\omega f= \tmap 0\omega \pistar \alpha\omega g = \tmap 0\alpha \tmap \alpha\omega \pistar \alpha\omega g = \tmap 0\alpha g,
\]
and $\tmap 0\alpha$ is in $I$ because $g$ is in $\sf_\alpha\subset \ia=T^{-1}_{0,\alpha}(I)$.  Thus $f$ is in $\ti$, as desired.

To prove the density of $\sf$ in $\ti$, first note that $\pistar \alpha\omega(\sf_\alpha)$ is dense in $\pistar \alpha\omega(\ia)$, so it suffices to show that $\bigcup_{0\leq\alpha<\omega} \pistar \alpha\omega(\ia)$ is dense in $\ti$.  Fix $f\in\ti$ and $\vep>0$ arbitrary.  We will show that $\|\pistar \alpha\omega(\tmap \alpha\omega f)-f\|<\vep$ for some $\alpha<\omega$.  Since $\tmap 0\alpha(\tmap \alpha\omega f) = \tmap 0\omega f$ is in $I$, the function $\tmap\alpha\omega f$ is in $\ia$, so this will establish the desired density.

By the definition of $\ta=A_\omega$, there exists $\alpha<\omega$ and $k\in \aa$ such that
\begin{equation}\label{first-numbered-eq}
\|f-\pistar \alpha\omega k\|<\vep/2.
\end{equation}
Then
\[
\|(\pistar\alpha\omega\circ \tmap \alpha\omega)(f) - (\pistar\alpha\omega\circ \tmap \alpha\omega)(\pistar \alpha\omega k) \|<\vep/2.
\]
Since $\tmap\alpha\omega \circ \pistar\alpha\omega$ is the identity, this gives
\begin{equation}\label{second-numbered-eq}
\|(\pistar\alpha\omega\circ \tmap \alpha\omega)(f) - (\pistar\alpha\omega k) \|<\vep/2.
\end{equation}
From (\ref{first-numbered-eq}) and (\ref{second-numbered-eq}) we get
\[
\|(\pistar \alpha\omega \circ \tmap\alpha\omega)(f) - f\|<\vep.
\]
This concludes the proof of condition (iv$'$).

A simple compactness argument shows that for every point $\tilde x\in \tx\sm\ts$ and every compact subset $E$ of $\tx\sm\ts$, there exists a neighborhood $U$ of $x$ and a function $f\in\sf$ such that $f|U=1$ and $f|E=0$.  Consequently, condition (v) can be proven in the same manner as was done earlier when we took $\sf_\alpha=\ia$.
\epf

%
%
%
%

\section{Normal peak point algebras not strongly regular}\label{peak-point-algebras}

In this section we prove Theorems~\ref{peak-point} and~\ref{peak-point-mod}.

\bpf[Proof of Theorem~\ref{peak-point}]
Set $A=R(K)$ with $K$ as in Theorem~\ref{main-theorem}, and set $I=\ol J_1$.  Then let $B$ be the uniform algebra $\ta$ obtained by applying Theorem~\ref{Cole-Izzo-thm2} taking $\tx$ to be metrizable.    
Let $x_0$ be the unique point of $\tpi^{-1}(1)$.  The uniform algebra $B$ is normal, $B$ has bounded relative units at every point of $\tx\sm\{x_0\}$, and every point of $\tx\sm\{x_0\}$ is a peak point for $B$.  The point $x_0$ is a peak point for $B$ also because the function 
$\bigl((1+z)/2\bigr)\circ\tpi$
peaks at $x_0$.

There is an order-preserving bijection between the closed ideals of $B$ containing $\ti$ and the closed ideals of $R(K)$ containing $I$.  Thus the family of ideals $\{I_\rho: 0<\rho<\infty\}$ in $R(K)$ yields the family of ideals $\{H_\rho:0<\rho<\infty\}$ in $B$.
\epf

\bpf[Proof of Theorem~\ref{peak-point-mod}]
The proof is essentially the same as the previous proof except that now we set $A=R(K)$ with $K$ as in Theorem~\ref{main-theorem-mod}, let $\Lambda$ be as in Theorem~\ref{main-theorem-mod}, and set $I=\ol {J_\Lambda}$.
\epf

\end{document}